\documentclass[preprint]{elsarticle}
\biboptions{numbers,sort&compress}
\usepackage[T1]{fontenc}
\usepackage{lmodern}
\usepackage{amsmath,amssymb,amsthm,mathtools,mathrsfs}
\usepackage[top=25mm,bottom=25mm,left=30mm,right=25mm]{geometry}
\usepackage[hidelinks]{hyperref}
\usepackage{xcolor}

\allowdisplaybreaks[3]
\numberwithin{equation}{section}
\mathtoolsset{showonlyrefs}
\newcommand{\R}{\mathbb R}
\newcommand{\s}{\mathbb S}
\newcommand{\N}{\mathbb N}
\newcommand{\abs}[1]{\lvert#1\rvert}
\newcommand{\norm}[1]{\lVert#1\rVert}
\newcommand{\LL}{\mathscr L}
\newtheorem{theorem}{Theorem}[section]
\newtheorem{introtheorem}{Theorem}

\newtheorem{introtheoremB}{Theorem}

\newtheorem{lemma}[theorem]{Lemma}
\newtheorem{proposition}[theorem]{Proposition}

\newtheorem{definition}[theorem]{Definition}
\theoremstyle{remark}
\newtheorem{remark}[theorem]{Remark}
\title{Nondegeneracy of low-energy nodal solutions to the Yamabe equation}

\author{Sérgio L. N. Neves}
\address{Department of Mathematics, Institute of Biosciences, Letters and
Exact Sciences (IBILCE), S\~ao Paulo State University (UNESP), Rua Crist\'ov\~ao
Colombo 2265, Jardim Nazareth, 15054-000 S\~ao Jos\'e do Rio Preto, SP, Brazil}
\author{Angela Pistoia}
\address{Department of Basic and Applied Sciences for Engineering,
Sapienza University of Rome, Via Antonio Scarpa 16, 00161 Rome, Italy}
\date{}
\begin{document}
\begin{abstract}
We study the linearized Yamabe equation in $\R^n$, or equivalently on
$\s^n$, at a low-energy sign-changing solution with two nodal domains.
We consider odd dimensions $n=2k-1$ and a solution invariant under
$\Gamma=O(k)\times O(k)$. Its regular level sets are isoparametric
hypersurfaces, which foliate the sphere away from the two focal
submanifolds. This symmetry reduces the equation for the solution to
a singular ODE. We choose a single-node solution with a profile
that has an extra symmetry to prove nondegeneracy in the
$\Gamma$-invariant subspace for every $k\ge2$, and in the full space
$H^1(\s^n)$ for sufficiently large $k$. Consequently, for such $k$, every kernel element is generated
by a conformal invariance of the equation. The proof uses separation
of variables, Wronskian identities and suitable comparison functions.
\end{abstract}
\maketitle

\section{Introduction}

We consider the critical elliptic equation
\begin{equation} \label{eqRn}
 -\Delta u=\frac{n(n-2)}4|u|^{\frac4{n-2}}u
 \quad\text{in }\R^n,\qquad
 u\in D^{1,2}(\R^n),\qquad n\ge 3,
\end{equation}
where $D^{1,2}(\R^n)$ is the completion of $C_c^\infty(\R^n)$ with
respect to the norm $\norm{u}=\sqrt{\int_{\R^n}\abs{\nabla u}^2}$. 

It is well known (see \cite{aubin,cgs,obata,talenti}) that all positive solutions of \eqref{eqRn}  are given by
the so-called {\em bubble}
\begin{equation} \label{famRn}
 U(x)
 =\left(\frac{2}{1+|x|^2}\right)^{\frac{n-2}{2}}
\end{equation}
and its translations and dilations
$$U_{\delta,\xi}(x)=\delta^{-\frac{n-2}2}U\left(\frac{x-\xi}\delta\right), \quad \delta>0,\quad \xi\in\R^n.$$

For sign-changing solutions, problem \eqref{eqRn} exhibits a much greater variety of profiles and nodal configurations.  The conformal equivalence with the round sphere $(\s^n,g)$ provides a useful way
to construct nodal solutions and study their geometry. Under
stereographic projection, the equation \eqref{eqRn} corresponds to the Yamabe equation on the sphere
\begin{equation} \label{eqSn}
 -\Delta_g v+\frac{n(n-2)}4v
 =\frac{n(n-2)}4|v|^{p-1}v,
 \qquad v \in H^1(\s^n),\qquad p=\frac{n+2}{n-2}.
\end{equation}

The first existence result for sign-changing solutions is due to Ding~\cite{Ding}, who proved, by variational methods, the existence of infinitely many nodal solutions to \eqref{eqSn} invariant under the action of the group $\Gamma = O(k)\times O(n+1-k) \subset O(n+1)$, $2\le k\le n-1$. Subsequently,
Clapp~\cite{clapp} proved the existence of a least-energy solution in a suitable equivariant class (different from the solutions found by Ding when $n\ge4$), whose energy was estimated by Clapp, Pistoia and Weth~\cite{clapiwe}.

 A more explicit construction was
later given by del Pino, Musso, Pacard and Pistoia~\cite{delpimupapi2011,delpimupapi2013}, who used
finite-dimensional reduction to build families of sign-changing solutions with large energy that concentrate along particular submanifolds of $\s^n$. For example, if $n\ge4$, they obtain sequences of solutions whose energy concentrates along one great circle or finitely many linked great circles, if $n\ge5$, they construct sequences of solutions whose energy concentrates along a two-dimensional torus. Other types of solutions have been constructed by Medina, Musso and Wei in \cite{mmw}, by Medina and Musso in \cite{mm} and, more recently, by Li and Sun in \cite{LiSun}.

Fern\'andez and Petean~\cite[Theorem~1.2]{FP}  constructed solutions to \eqref{eqSn}
whose nodal sets consist of a prescribed number of {\em isoparametric
hypersurfaces.} More precisely, for an isoparametric function on $\s^n$ whose regular level sets have at least two distinct principal curvatures, they reduce \eqref{eqSn} to a singular ODE on $[0,\pi]$ with Neumann boundary conditions. For any positive integer $\kappa$, they obtain a profile with exactly $\kappa$ zeros, yielding a solution of \eqref{eqSn} whose nodal set consists of exactly $\kappa$ connected isoparametric hypersurfaces. For $\kappa=1$, they obtain a \emph{single-node} solution (see Definition~\ref{onu2}).
Every nodal solution of \eqref{eqRn} has energy strictly larger than twice the energy of the bubble, see Weth~\cite{Weth}, since each nodal domain carries at least the energy of $U$. Solutions with exactly two nodal domains are therefore the natural candidates for nodal solutions of lowest energy, and in this sense the single-node solutions are \emph{low-energy} nodal solutions, in the terminology of~\cite{FP}. We refer to~\cite{FP} for energy estimates of these solutions, and to~\cite{clapiwe} for an upper bound on the least energy of a nodal solution of \eqref{eqSn}.

The purpose of this paper is to investigate the nondegeneracy of these
single-node solutions in the sense introduced by
Duyckaerts, Kenig and Merle~\cite{DKM}. We recall this notion through the
kernel elements generated by the geometric invariances of the equation.
Let $u\in D^{1,2}(\R^n)$ be a nonzero solution of \eqref{eqRn}, and define
\[
 L_u=-\Delta-\frac{n(n+2)}4|u|^{\frac4{n-2}},
 \qquad
 \mathcal Z_u=\{\psi\in D^{1,2}(\R^n):L_u\psi=0\}.
\]
The equation \eqref{eqRn} is invariant under translations, dilations, orthogonal
transformations and the Kelvin transform. Differentiating the equation along these invariances, we obtain the following
 explicit elements of $\mathcal Z_u$ (see \cite[p.~734, formulas~(1.6)--(1.8)]{DKM}):
\[
 \begin{aligned}
 D_u(x)&=\frac{n-2}{2}u(x)+x\cdot\nabla u(x),\\
 T_{u,i}(x)&=\partial_{x_i}u(x),
                  &&1\le i\le n,\\
 R_{u,ij}(x)&=x_i\partial_{x_j}u(x)-x_j\partial_{x_i}u(x),
                  &&1\le i<j\le n,\\
 C_{u,i}(x)&=(2-n)x_i u(x)+|x|^2\partial_{x_i}u(x)
                         -2x_i\,x\cdot\nabla u(x),
                  &&1\le i\le n.
 \end{aligned}
\]
If
\[
 \widetilde{\mathcal Z}_u
 =\operatorname{span}\bigl(\{D_u\}
 \cup\{T_{u,i},C_{u,i}:1\le i\le n\}
 \cup\{R_{u,ij}:1\le i<j\le n\}\bigr),
\]
 it is clear that $\widetilde{\mathcal Z}_u\subseteq\mathcal Z_u$.
The nondegeneracy condition stated by Duyckaerts, Kenig and Merle~\cite{DKM} is precisely
\[
 \mathcal Z_u=\widetilde{\mathcal Z}_u,
\]
that is, every finite-energy solution of the linearized equation is a
linear combination of the variations induced by these symmetries.
There are
\[
 N=2n+1+\frac{n(n-1)}2=\frac{(n+1)(n+2)}2
\]
generators in the above list, and hence
$\dim\widetilde{\mathcal Z}_u\le N$.
This dimension can be strictly smaller when the solution has
continuous symmetries, since some infinitesimal variations then vanish
or are linearly dependent. Thus $N$ is the maximal dimension allowed
by the conformal invariances for a geometrically nondegenerate
solution and nondegeneracy does not require that this maximum be attained.

The positive bubbles provide a basic example of nondegeneracy in this
sense. For $u=U_{\delta,\xi}$, the finite-energy solutions of
\begin{equation}\label{lineqRn}
 -\Delta\psi=\frac{n(n+2)}4u^{\frac4{n-2}}\psi
 \quad\text{in }\R^n
\end{equation}
are precisely the linear combinations of the derivatives with respect
to the parameters in \eqref{famRn}, see \cite[Lemma~A.1]{BE}. Hence
\[
 \mathcal Z_u=\operatorname{span}\{\partial_\delta U_{\delta,\xi},
       \partial_{\xi_1}U_{\delta,\xi},\ldots,
       \partial_{\xi_n}U_{\delta,\xi}\}
       =\widetilde{\mathcal Z}_u,
 \qquad \dim\mathcal Z_u=n+1.
\]
Thus the
positive bubble is nondegenerate in the sense of
Duyckaerts--Kenig--Merle, although its kernel dimension is smaller
than $N$.


The nondegeneracy of sign-changing solutions was first studied by Musso and Wei~\cite{MW}, who proved nondegeneracy results for the {\em crown solutions} constructed by del Pino, Musso, Pacard and Pistoia~\cite{delpimupapi2011}. Medina, Musso and Wei~\cite{mmw} established nondegeneracy for the solutions they constructed. These arguments exploit the multi-bubble structure arising from the Lyapunov--Schmidt reduction, combining the nondegeneracy of the positive bubble with the symmetries of the configuration. Medina and Musso~\cite{mm} constructed a different family and proved that its geometric symmetry directions have maximal rank in dimension three, nondegeneracy of that family remains open in their work. More recently, Li and Sun~\cite{LiSun} constructed two families of nodal solutions concentrating along two planar circles, one of which is conformally equivalent to the family of Medina and Musso, and showed that these solutions also have maximal rank in dimension three.

In the present paper, we establish nondegeneracy, in all sufficiently large odd dimensions, for the \emph{symmetric single-node} solution of problem \eqref{eqRn} introduced in Definition~\ref{onu2} (see also Remark~\ref{onu3}). Its profile is described below and can be obtained from the construction of Fern\'andez and Petean~\cite{FP}, as explained in Section~\ref{sec:setting}.
Let us state our main result.
\begin{introtheorem}\label{thm:euclidean}
There exists $k_0\ge2$ such that, for every $k\ge k_0$, the following holds.
Let $n=2k-1$ and let $u$ be the symmetric single-node solution of problem \eqref{eqRn}.
Then the following $k(k+2)$ functions form a basis of the linearized kernel in $D^{1,2}(\R^n)$ at $u$:
\begin{equation}\label{eq:euclidean-basis}
 \begin{aligned}
 Z_0(z)&=\frac{n-2}{2}u(z)+z\cdot\nabla u(z),\\
 Z_i(z)&=\partial_{z_i}u(z),
       &&1\le i\le n,\\
 Z_{n+i}(z)&=-2z_iZ_0(z)+|z|^2Z_i(z),
       &&1\le i\le k,\\
 Z_{n+k+(i-1)(k-1)+a}(z)
       &=-z_{k+a}Z_i(z)+z_iZ_{k+a}(z),
       && {1\le i\le k,\ 1\le a\le k-1.}
 \end{aligned}
\end{equation}
In particular, with
$L_u=-\Delta-\frac{n(n+2)}4|u|^{4/(n-2)}$,
\[
 \ker_{D^{1,2}(\R^n)}L_u
 =\operatorname{span}\{Z_\alpha:0\le\alpha\le k(k+2)-1\}
 =\widetilde{\mathcal Z}_u,
\]
\[
 \dim\ker L_u=k(k+2)=\frac{(n+1)(n+5)}4.
\]
Thus $u$ is nondegenerate in the sense of Duyckaerts--Kenig--Merle.
\end{introtheorem}

To introduce the {\em single-node solution},
let us briefly recall the construction of Fern\'andez and Petean, in which an isoparametric function on the sphere is used to reduce the equation \eqref{eqSn} to an ODE.
The theory of isoparametric hypersurfaces is very rich and dates back to E. Cartan in the late 1930s with important contributions by H. M\"unzner in the early 1980s, the recent monograph \cite{CR} is an excellent reference on the subject. To recall the basic definitions, a smooth function $F : \s^n \to \R$ is called an {\em isoparametric function} if $\abs{\nabla_g F}^2$ and $\Delta_g F$ are smooth functions of $F$ itself, so $\abs{\nabla_g F}^2 = \mathcal A \circ F$ and $\Delta_g F = \mathcal B \circ F$ for some smooth real functions $\mathcal A$ and $\mathcal B$. The regular level sets form an isoparametric family of hypersurfaces, and the singular levels are the focal submanifolds. An oriented hypersurface in the sphere is isoparametric precisely when its principal curvatures are constant. 

In this case, if $W$ is a smooth real function of one variable such that $W(F(y))$ is a smooth function of $y \in \s^n$, we have the identity
\begin{equation}
 \Delta_g (W \circ F) = (W'' \circ F)\abs{\nabla_g F}^2 + (W' \circ F)\Delta_g F   
\end{equation}
and problem \eqref{eqSn} reduces to the ODE
\begin{equation}\label{odeW}
  \mathcal A(x)W'' + \mathcal B(x)W' + \frac{n(n-2)}{4} \left(\abs{W}^{p-1}W - W\right) = 0 \qquad x=F(y).
\end{equation}

Now consider 
\[
 \R^{n+1}=\R^k\times\R^{n+1-k},\qquad 2\le k\le n-1,
\]
and the Cartan--M\"unzner polynomial $F:\R^{n+1}\to\R$ defined by 
\begin{equation*}
  F(y',y'') = \abs{y'}^2 - \abs{y''}^2 \qquad y=(y',y '') \in \R^k \times \R^{n+1-k} = \R^{n+1}. 
\end{equation*} 
The restriction $F |_{\s^n}$ is an isoparametric function (see \cite{CR}) and corresponds to the case where each regular level set has two distinct principal curvatures with multiplicities $m_1 = k-1$ and $m_2=n-k$. On $\s^n$ it follows that
\[
 |\nabla_g F|^2=4(1-F^2),\qquad \Delta_g F=2(2k-n-1)-2(n+1)F,
\]
and then
\begin{equation*}
 4(1-x^2)W''+(2(2k-n-1)-2(n+1)x)W'+\dfrac{n(n-2)}{4}(|W|^{p-1}W-W)=0,\qquad -1<x<1.
\end{equation*}
Performing the further change of variables $w(r)=W(\cos r)$, one sees that, if
$w$ solves the boundary value problem
\begin{equation}\label{bvp}
 \begin{cases}
 w''+\dfrac{h(r)}{\sin r}w'
 +\dfrac{n(n-2)}{16}\bigl(|w|^{p-1}w-w\bigr)=0,
       &r\in(0,\pi),\\[2mm]
 w'(0)=w'(\pi)=0,
 \end{cases}
\end{equation}
where
\begin{equation} \label{coefh}
 h(r)=\frac{n-1}{2}\cos r-\frac{2k-n-1}{2},
\end{equation}
then the function
\begin{equation}\label{vonu}
 v(y)=w\left( \arccos(|y'|^2-|y''|^2)\right) \qquad y=(y',y'')\in\s^n\subset\R^k\times\R^{n+1-k}
\end{equation}
solves the Yamabe equation \eqref{eqSn}. Notice that the boundary conditions in \eqref{bvp} are necessary for regular bounded solutions.
Moreover, if $u$ denotes its stereographic image with pole $e_{n+1}$, namely
\begin{equation}\label{onu}
 u(z)=\left(\frac{2}{1+|z|^2}\right)^{\frac{n-2}{2}}
 w\left(\arccos\left(\frac{8|z'|^2}{(1+|z|^2)^2}-1\right)\right),\quad z=(z',z'')\in\R^k\times\R^{n-k},
\end{equation}
then $u$ is a finite-energy solution of \eqref{eqRn}.

By \cite[Theorem~1.2]{FP}, for any integer $\kappa\ge1$ problem \eqref{bvp} admits a sign-changing solution $w_\kappa$ having exactly $\kappa$ zeros in $(0,\pi)$.

With a slight abuse of terminology, we shall refer to a solution $v$ to \eqref{eqSn}, defined in \eqref{vonu}, as \emph{radial}.
Indeed, the isoparametric foliation under consideration
consists of tubes around either of its two focal submanifolds.
Hence, in this case the function $v$ depends only on the
geodesic distance to a fixed focal submanifold, extending
the usual notion of radial symmetry about a point.

\begin{definition}\label{onu2}
We call the function $u$ defined in \eqref{onu} a \emph{single-node solution} of \eqref{eqRn} if its profile $w$ has exactly one zero in $(0,\pi)$. The corresponding solutions in $\R^n$ or $\s^n$ have exactly two nodal domains.
If, in addition, $n=2k-1$ is odd and the profile satisfies
\begin{equation} \label{eq:solution}
 w(\pi-r)=-w(r),\qquad
 w>0\ \text{on }(0,\pi/2),\qquad
 w'<0\ \text{on }(0,\pi),
\end{equation}
we call $u$ a \emph{symmetric single-node solution} of \eqref{eqRn}.
\end{definition}

\begin{remark}
It would also be interesting to compare the single-node solutions considered here with the solutions obtained by Ding~\cite{Ding} through variational methods in the same invariant class. In particular, does the single-node profile occur among those variational solutions?
\end{remark}

\begin{remark}\label{onu3} If $n=2k-1$ is odd, an extra symmetry is obtained by performing the construction with $\R^{n+1}=\R^k\times\R^k$. Then the constant term in $h$ defined in \eqref{coefh} vanishes and equation \eqref{bvp} is invariant under $w(r)\mapsto-w(\pi-r)$. In Section~\ref{sec:setting} we explain how to choose a single-node profile $w$ whose unique zero is at $\pi/2$ and which satisfies the conditions of Definition~\ref{onu2}, in particular, a symmetric single-node solution exists and, by Lemma~\ref{lem:unique-profile}, it is unique.
The corresponding spherical solution $v$ defined in \eqref{vonu} has exactly one nodal hypersurface, equal to the product of spheres
\[
 \s^{k-1}(1/\sqrt2)\times\s^{k-1}(1/\sqrt2),
\]
and the nodal set of $u$ is its stereographic image,
\[
 \{z\in\R^{2k-1}:(1+|z|^2)^2=8|z'|^2\}.
\]
\end{remark}

We now provide some comments on our main result and its proof. First, we remark that the dimension of our kernel is strictly smaller than the maximal
conformal dimension $N$. More precisely,
\[
 N-\dim\ker L_u
 =\frac{(n+1)(n+2)}2-\frac{(n+1)(n+5)}4
 =\frac{n^2-1}4=k(k-1)>0.
\]
This deficit has a direct geometric interpretation, the subgroup
$O(k)\times O(k)$ fixing the solution on the sphere has dimension
$k(k-1)$, and its infinitesimal rotations produce no variation of
the solution.
The elements of the basis \eqref{eq:euclidean-basis} correspond to one
dilation, $n$ translations, $k$ special conformal directions, and
$k(k-1)$ rotations mixing the two coordinate blocks in the decomposition $\R^{n}=\R^k\times\R^{k-1}$.

In our setting, the single-node solution is not obtained as a perturbation
of a configuration of standard bubbles, and the analysis cannot be
reduced to the well-understood linearized operator at the positive
bubble as in the papers \cite{MW,mmw}. Instead, we must study the linearized equation directly
around a non-explicit sign-changing profile. This requires new
arguments for the singular ordinary differential equations
arising from the isoparametric decomposition, in order to exclude
additional kernel elements beyond those generated by conformal
invariance. We establish this result in sufficiently large odd dimensions.

The proof uses separation of variables and spherical harmonics
on the two factors $\s^{k-1}$ to decompose the linearized equation into a system of ordinary differential equations. For every $k\ge2$, we prove radial nondegeneracy, i.e., nondegeneracy restricted to the subspace of functions invariant under the action of the subgroup $O(k) \times O(k)$. To prove general nondegeneracy, we need the restriction to sufficiently large $k$.  {More precisely, the following results hold in every odd dimension.}

\begin{introtheoremB}\label{thm:B}
Let $k\ge2$, $n=2k-1$ and $\Gamma=O(k)\times O(k)$. Then:
\begin{enumerate}
\item[(i)] there is exactly one solution $w$ of \eqref{eq:w} satisfying \eqref{eq:solution}; in particular, the symmetric single-node solution is unique;
\item[(ii)] the corresponding solution $v$ of \eqref{eqSn} is nondegenerate in the $\Gamma$-invariant subspace, namely $\ker L_v\cap H^1_\Gamma(\s^n)=\{0\}$.
\end{enumerate}
\end{introtheoremB}


\begin{remark} Extending the result to a different decomposition $\R^{n+1}=\R^k\times\R^{n+1-k}$ presents two main difficulties. The first already appears in the radial mode. When $k=(n+1)/2$, the profile $w$ chosen here has its unique zero at $\pi/2$, which is also the zero of $h$.
Moreover, $W(x)=w(\arccos x)$ is odd and $\abs{W}$ is convex. These properties
are heavily used throughout the argument. For an unequal decomposition,
the reflection symmetry is lost, the zero of $h$ is displaced, and the
nodal point is no longer fixed by symmetry. 

A further difficulty arises in the angular modes. Separation into
spherical harmonics is still available for general $k$, but the two
factors contribute the eigenvalues
\[
 j(j+k-2)\qquad\text{and}\qquad
 \ell(\ell+n-k-1),
\]
respectively. In the symmetric case $k=(n+1)/2$, reflection interchanges the two
modes $(j,0)$ and $(0,\ell)$, allowing them to be treated by the same
argument. When the dimensions differ, the estimates obtained here
do not directly control both cases, an extension would
require new arguments adapted to the asymmetric case.
\end{remark}

The paper is organized as follows. In Section~\ref{sec:setting} we specify the solution and state the equivalent kernel classification on the sphere. In Section~\ref{sec:separation} we use spherical harmonics and separation of variables to derive the system of equations for the functions in the kernel. The subsequent sections are devoted to the proofs of radial nondegeneracy in the
$\Gamma$-invariant subspace $H^1_\Gamma(\s^{2k-1})$,  with $\Gamma=O(k)\times O(k)$, for every $k\ge2$, and of nondegeneracy in the full space $H^1(\s^{2k-1})$ for sufficiently large $k$. Finally, Section~\ref{sec:euclidean-kernel} completes the proofs of Theorem~\ref{thm:euclidean} and Theorem~\ref{thm:B}.

\section{The Yamabe equation and the single-node solution}\label{sec:setting}
The exponent and coefficients associated with the fixed parameter $k$ are
\begin{equation}\label{eq:constants}
 \begin{gathered}
 n=2k-1,\qquad p=\frac{2k+1}{2k-3},\qquad
 b=p-1=\frac4{2k-3},\\
 c=\frac{(2k-1)(2k-3)}{16},\qquad
 d=bc=\frac{2k-1}{4}.
 \end{gathered}
\end{equation}
Every solution, function, coefficient,
operator, and constant in this paper depends on this fixed $k$. To keep the
notation simple, we suppress this dependence. We also set
\[
 f(s)=c\left(|s|^{p-1}s-s\right),\qquad s\in\R.
\]

On the round sphere $(\s^n,g)$ we consider the Yamabe equation
\begin{equation}\label{eq:yamabe}
 -\Delta_g v+4cv=4c|v|^{p-1}v.
\end{equation}
Write $y=(y',y'')\in\R^k\times\R^k$ and set
$F(y)=|y'|^2-|y''|^2$. On $\s^{2k-1}$,
\[
 |\nabla_g F|^2=4(1-F^2),\qquad \Delta_g F=-4kF.
\]
Set
\[
 x:=F(y)=|y'|^2-|y''|^2\in[-1,1],
 \quad r:=\arccos F(y),\quad
\hbox{i.e.} \ x=\cos r. \]
 The ansatz $v(y)=W(x)=W(F(y))$
therefore reduces \eqref{eq:yamabe} to
\begin{equation}\label{eq:W}
 (1-x^2)W''-kxW'+c(|W|^{p-1}W-W)=0,\qquad -1<x<1.
\end{equation}
Equivalently, setting $w(r)=W(\cos r)$,
\begin{equation}\label{eq:w}
 w''+(k-1)\cot r\,w'+c(|w|^{p-1}w-w)=0,
 \qquad w'(0)=w'(\pi)=0.
\end{equation}
From now on we fix a single-node solution $w$ with the properties \eqref{eq:solution}, its existence can be obtained from \cite{FP}. For completeness, let us sketch the argument. First, for $\beta>0$ let $w(\cdot;\beta)$ be the solution of the initial value problem
\begin{equation} \label{ivp}
\begin{cases}
w''+(k-1)\cot r\,w'+f(w)=0 \\
w(0)=\beta \\
w'(0)=0
\end{cases}
\end{equation}
By the results of \cite[Section~4]{FP} there exists $\beta_\ast > 1$ for
which the first zero of $w(\cdot;\beta_\ast)$ occurs at $\pi/2$. Odd reflection across $\pi/2$
then gives a solution on $[0,\pi]$ with $w'(\pi)=0$.
To check strict monotonicity, consider the energy function
\begin{equation} \label{energy}
 E(r)=\frac12w'(r)^2+c\left(\frac{|w(r)|^{p+1}}{p+1}
                              -\frac{w(r)^2}{2}\right).
\end{equation}
On $(0,\pi/2)$,
$E'(r)=-(k-1)\cot r\,w'(r)^2\le0$, and
$E(\pi/2)=w'(\pi/2)^2/2>0$ by uniqueness for the ODE, hence $E(r)\ge E(\pi/2)>0$ on $(0,\pi/2]$.
At any critical point of $w$ in this interval, the positivity of $E$ implies
$w^{p-1}>(p+1)/2>1$, hence $w''=-c(w^p-w)<0$.
Since $w'<0$ near $0$ (because $f(\beta_\ast)>0$), $w'$ cannot have a first zero in $(0,\pi/2)$.
Odd reflection proves $w'<0$ on the whole interval $(0,\pi)$.


For $v$ as in \eqref{vonu}, define the linearized operator of \eqref{eqSn} by
\[
 L_v\psi=-\Delta_g\psi+4c(1-p|v|^{b})\psi.
\]
Our main result is the following.
\begin{theorem}\label{thm:main}
There exists $k_0\ge2$ such that, for every $k\ge k_0$,
\begin{equation}\label{eq:kernel}
 \ker_{H^1(\s^{2k-1})}L_v
 =\operatorname{span}\{
 h_{10}(r)\theta_i,\ h_{01}(r)\omega_a,
 \ w'(r)\theta_i\omega_a:1\le i,a\le k\},
\end{equation}
where $y=(\cos(r/2)\theta,\sin(r/2)\omega)$,
$\theta,\omega\in\s^{k-1}$, and
\begin{equation}\label{eq:geometric-modes}
 \begin{aligned}
 h_{10}(r)&=\frac{2k-3}{2}w(r)\cos\frac r2
                  +2w'(r)\sin\frac r2,\\
 h_{01}(r)&=\frac{2k-3}{2}w(r)\sin\frac r2
                  -2w'(r)\cos\frac r2.
 \end{aligned}
\end{equation}
In particular, $\dim\ker L_v=k(k+2)$.
\end{theorem}
The proof is completed in Section~\ref{sec:one-sided}.
Theorem~\ref{thm:euclidean} identifies this kernel with the infinitesimal
conformal invariances of the corresponding Euclidean solution.

\section{Separation of variables and convexity}\label{sec:separation}
In the coordinates of Theorem~\ref{thm:main},
\[
 \Delta_g=4\left(\partial_r^2+(k-1)\cot r\,\partial_r\right)
 +\frac{\Delta_{\s^{k-1},\theta}}{\cos^2(r/2)}
 +\frac{\Delta_{\s^{k-1},\omega}}{\sin^2(r/2)}.
\]
Let $\mathcal H_j$ be the space of spherical harmonics of degree $j$ on $\s^{k-1}$,
with eigenvalue
\[
 \lambda_j=j(j+k-2),\qquad j=0,1,\ldots,
\]
for $-\Delta_{\s^{k-1}}$, and let $\mathfrak d_j=\dim\mathcal H_j$.
Choosing orthonormal bases $\{H_{j,m}\}_{m=1}^{\mathfrak d_j}$ of $\mathcal H_j$ (so that $\{H_{j,m}\}_{j,m}$ is an orthonormal basis of $L^2(\s^{k-1})$), we expand
\[
 \psi=\sum_{j,\ell\ge0}\sum_{m=1}^{\mathfrak d_j}\sum_{m'=1}^{\mathfrak d_\ell}
 h_{j,\ell,m,m'}(r)H_{j,m}(\theta)H_{\ell,m'}(\omega).
\]
If $L_v\psi=0$, then every radial coefficient $h_{j,\ell,m,m'}$ solves
\begin{equation}\label{eq:mode}
 h''+(k-1)\cot r\,h'+f'(w)h
 -\left(\frac{\lambda_j}{4\cos^2(r/2)}
       +\frac{\lambda_\ell}{4\sin^2(r/2)}\right)h=0.
\end{equation}
Indeed, since $|v|^{b}$ is H\"older continuous, elliptic regularity gives $\psi\in C^{2}(\s^n)$. The series converges in $L^2(\s^n)$, and each coefficient
\[
 h_{j,\ell,m,m'}(r)=\int_{\s^{k-1}}\int_{\s^{k-1}}\psi\,H_{j,m}(\theta)H_{\ell,m'}(\omega)\,d\theta\,d\omega
\]
is bounded on $[0,\pi]$ and of class $C^2$ on $(0,\pi)$; equation \eqref{eq:mode} follows by testing $L_v\psi=0$ against $H_{j,m}(\theta)H_{\ell,m'}(\omega)$ and using the symmetry of $\Delta_{\s^{k-1}}$. Moreover, $\psi\equiv0$ if and only if all the coefficients vanish. The indicial roots of \eqref{eq:mode} are $\ell$ and $-(\ell+k-2)$ at $r=0$, and $j$ and $-(j+k-2)$ at $r=\pi$, when $k=2$ and $\ell=0$ (respectively $j=0$) the root $0$ is double and the second solution has a logarithmic singularity. Hence boundedness selects, at each endpoint, a one-dimensional space of solutions, which we call \emph{regular}.
Regular solutions satisfy $h(r)=O(r^\ell)$ near $0$ and
$h(r)=O((\pi-r)^j)$ near $\pi$, if $\ell=0$ then $h'(0)=0$ and likewise $h'(\pi)=0$ if $j=0$.  
Here $f$ is the nonlinearity introduced in Section~\ref{sec:setting}.

We fix the notation
\begin{equation}\label{eq:operators}
 a(r)=\sin^{k-1}r,\qquad
 \LL h=-(ah')'-a f'(w)h,\qquad
 V(x)=c(p|W(x)|^{b}-1).
\end{equation}
Then \eqref{eq:mode} is
\[
 \LL h+ a\left(\frac{\lambda_j}{4\cos^2(r/2)}
              +\frac{\lambda_\ell}{4\sin^2(r/2)}\right)h=0.
\]
Throughout, $V$ is always regarded as a function of the variable $x$. First, we prove some properties of the transformed solution $W$.

\begin{lemma}\label{lem:convex}
For every $k\ge2$, $W$ is odd and
\[
 Y:=W'>0\quad\text{on }[-1,1],\qquad
 W''>0,\quad W<xW'\quad\text{on }(0,1).
\]
Moreover, $Y$ is even and satisfies
\begin{equation}\label{eq:Y}
 (1-x^2)Y''-(k+2)xY'+(V-k)Y=0.
\end{equation}
\end{lemma}
\begin{proof}
Oddness and positivity of $W'$ follow from \eqref{eq:solution}. At $x=1$,
\eqref{eq:W} gives $W'(1)=c(W(1)^p-W(1))/k>0$.
Differentiation gives \eqref{eq:Y}. On $(0,1)$ set
$J=(1-x^2)W''$ and $\vartheta=(1-x^2)^{k/2}$. Then
\[
 (\vartheta J)'=\vartheta W'(k+c-cpW^{b}).
\]
The factor in parentheses strictly decreases and is positive at $0$.
Also $(\vartheta J)(0)=0$ and $\vartheta J\to0$ at $1$, by regularity
at the pole. If it remained positive on $(0,1)$, then $\vartheta J$ would be strictly
increasing, which is incompatible with $(\vartheta J)(0)=0$ and
$\vartheta J\to0$ at $1$. Hence this factor changes sign exactly once, thus
$\vartheta J$ first increases and then decreases to zero. Thus $J>0$.
Finally,
$W(x)=\int_0^xW'(s)\,ds<xW'(x)$.
\end{proof}

\section{The mode \texorpdfstring{$(0,0)$}{(0,0)}: radial nondegeneracy}\label{sec:radial}
Let $\varphi$ be a solution to \eqref{eq:mode} with $j=\ell=0$, i.e.
\begin{equation}\label{eq:moderad}
\varphi''+(k-1)\cot r\,\varphi'+f'(w)\varphi =0,\ \varphi'(0)=\varphi'(\pi)=0.
\end{equation}

We aim to prove that $\varphi\equiv0.$

Suppose, for a contradiction, that there is a nontrivial solution. We normalize it by $\varphi(0)=1$ and prove that
$\varphi'(\frac{\pi}{2})\not =0$ and $\varphi(\frac{\pi}{2})\not=0$. A contradiction will follow (see Proposition~\ref{prop:radial}).

The proof separates these two tasks. First we determine the zero pattern
of $\varphi$, more precisely, we show that it has exactly one zero in $[0,\frac{\pi}{2}]$, which is an interior point, in
particular, $\varphi(\frac{\pi}{2})<0$. A second comparison is needed to determine
the sign of $\varphi'(\frac{\pi}{2})$. The functions introduced next provide the sign information needed to build that comparison.

\subsection{Auxiliary functions}
The quantities introduced below are all functions of the fixed solution $w$ and the coefficients defined in \eqref{eq:constants}. Since $w>0$ and $w'<0$ on
$(0,\frac{\pi}{2})$, the expressions are well defined there. 


 On
$(0,\frac{\pi}{2})$ set
\begin{equation}\label{eq:phase-def}
 \begin{gathered}
 t(r)=\tan^2r,\qquad S(r)=\left(w(r)\right)^{b},\\
 \xi(r)=-b\tan r\,\frac{w'(r)}{w(r)},\qquad
 \mu(r)=-\tan r\,\frac{f(w(r))}{w'(r)},\\
 G(r)=b\mu(r)+\xi(r)-b(k-1).
 \end{gathered}
\end{equation}
 
Direct differentiation gives
\begin{align}
 \mu(r)\xi(r)&=d t(r)(S(r)-1),\label{eq:muxi}\\
 \tan r\,\xi'(r)&=\xi(r)\left(\mu(r)+\frac{\xi(r)}{b}+t(r)-k+2\right)
              =\xi(r)\left(1+t(r)+\frac {G(r)}{b}\right),\label{eq:xi}\\
 \tan r\,\mu'(r)&=\mu(r)\left(k+t(r)-\mu(r)-\frac p{b}\xi(r)\right)-dt(r).
 \label{eq:mu}
\end{align}
At the origin, $w(0)>1$, and one has
\[
 w'(r)=-\frac{f(w(0))}k r+O(r^3),\qquad
 \xi(r)=O(r^2),\quad \mu(r)\to k,\quad G(r)\to b\qquad\hbox{as}\ r\to0.
\]

\begin{lemma}\label{lem:phase}
For every $k\ge2$ and $0<r<\frac{\pi}{2}$, 
\[G(r)>0,\quad \xi'(r)>0,\quad\mu(r)<k.\]
\end{lemma}
\begin{proof}
At a zero of $G$, \eqref{eq:xi}--\eqref{eq:mu} yield
\[
 \tan r\,G'(r)=t(r)+\xi^2(r)-\beta_0\xi(r)+\beta_0,
 \qquad\beta_0=b(k-1)=2+\frac{b}{2}.
\]
Since $0<\beta_0\le4$, the right-hand side equals
\[
 \tan^2 r+(\xi(r)-\beta_0/2)^2+\beta_0(1-\beta_0/4)>0.
\]
This excludes a first zero starting from $G(0^+)=b>0$.
Monotonicity of $\xi$ follows from \eqref{eq:xi}.
For the last assertion, Lemma~\ref{lem:convex} gives
\[
 \frac{d}{dr}\log Y(\cos r)
 =\frac{w''}{w'}-\cot r=(\mu-k)\cot r<0.
\]
\end{proof}

\begin{lemma}\label{lem:phase-bound}
For $k\ge4$ and $0<r<\frac{\pi}{2}$,
\begin{equation}\label{eq:phase-bound}
 \xi(r)<\tan^2r+\frac52.
\end{equation}
\end{lemma}
\begin{proof}
First, $\xi(r)-t(r)\to0$ as $r\to0$. Moreover, as $r\to \frac\pi2$,
\begin{equation}\label{eq:mid-expansion}
 w(r)=-w'\left(\frac{\pi}{2}\right)\varrho+O\left(\varrho^3\right),\qquad w'(r)=w'\left(\frac{\pi}{2}\right)+O\left(\varrho^2\right),
\end{equation}
where $\varrho:=\frac{\pi}{2}-r.$
Thus $\frac{\xi(r)}{t(r)}\to b<1$, and $\xi(r)-t(r)\to-\infty$ as $r\to\frac{\pi}{2}$.
If its maximum $M$ were at least $5/2$, it would be attained at an
interior point. At this point, write $t(r)=\tan^2r$ and $\xi(r)=t(r)+M$.
Since $\tan r\,t'=2t(1+t)$, \eqref{eq:xi} gives
\[
 \frac G{b}=\frac{(1+t)(t-M)}{t+M}>0,\qquad t>M.
\]
On the other hand, \eqref{eq:muxi} and $S>0$ imply
$\mu+dt/(t+M)>0$. Substitution of the preceding expression for $G$
gives instead
\[
 \mu+\frac{dt}{t+M}=-\frac{\mathfrak p(M,t)}{4(t+M)},
\]
where
\[
 \mathfrak p(M,t)=(2k-7)t^2+((4k-2)M-6k+1)t
              +(2k-3)M^2-4(k-2)M.
\]
For $k\ge4$, $M\ge5/2$ and $t>M$, this polynomial is strictly increasing
in $t$, and
\[
\mathfrak p(M,t)>\mathfrak p(M,M)
 =M((8k-12)M-10k+9)\ge M(10k-21)>0.
\]
This contradiction proves \eqref{eq:phase-bound}.
\end{proof}

\subsection{Comparison functions for the radial equation}
Let $\varphi$ be the solution of
\begin{equation}\label{eq:shoot}
 \varphi''+(k-1)\cot r\,\varphi'+f'(w)\varphi=0,
 \qquad\varphi(0)=1,\quad\varphi'(0)=0.
\end{equation}
This is the derivative of the solution $w(\cdot;\beta)$ in \eqref{ivp} with respect to $\beta$, evaluated at $\beta=w(0)$. The use of a
variation with respect to the shooting parameter and its zero pattern to prove uniqueness and nondegeneracy is a classical argument, see for instance Kwong~\cite{Kwong} and also Tang~\cite{Tang}, from whom we borrow some ideas.
Consider the linear operator $\mathcal D$ and the weighted Wronskian $\mathcal W$ defined by
\begin{align}
 \mathcal D h &= h''+(k-1)\cot r\,h'+f'(w)h, \label{opD}\\
 \mathcal W[g,h] &= a(h' g- h g'), \label{WWronsk}
\end{align}
where $a$ was defined in \eqref{eq:operators}. By differentiation we get the identity
\begin{equation} \label{WWronskIdentity}
 \left(\mathcal W[g,h] \right)'=a \left( g\mathcal D[h]-h\mathcal D[g] \right).
\end{equation}

To simplify notation, we write $\mathcal W_w:=  \mathcal W[\varphi,w] =a(w'\varphi-w\varphi')$ for the weighted Wronskian of $\varphi$ and $w$. 
Since $ \mathcal D w = wf'(w)-f(w)=d|w|^{p-1}w$, on $(0,\frac{\pi}{2})$ the Wronskian identity \eqref{WWronskIdentity} gives
\begin{equation}\label{eq:Wwprime}
 \mathcal W_w'=daw^p\varphi.
\end{equation}
In particular, $\mathcal W_w>0$ near the origin, because $\mathcal W_w(0)=0$ and $\varphi(0)=1$.

Notice that the derivative $\mathcal W_w'=d a w^p\varphi$ has the same sign as $\varphi$, we shall show that $\varphi$ changes sign exactly once on $(0,\frac{\pi}{2})$, so $\mathcal W_w$ starts increasing from $0$ and then decreases. Showing that it remains
positive up to $\frac{\pi}{2}$ is the key point. If we can find another function $z$ such that the expression $\mathcal D z$ is particularly simple, then we can use its Wronskian with $\varphi$ and the identity \eqref{WWronskIdentity} to get more information about $\varphi$. A good strategy, already used in \cite{KwongLi}, is to try combinations of $w$ and $w'$.


\begin{lemma}\label{lem:radial-comparison}
For every $k\ge2$ there are functions $z,\Phi$ on $(0,\frac{\pi}{2})$ such that, with
$\mathcal W_{z}:=\mathcal W[z,\varphi] =a(z\varphi'-\varphi   z   ')$,
\begin{equation}\label{eq:comparison-properties}
 \mathcal W_{z}'=\Phi\mathcal W_w',\qquad \Phi'<0,\qquad (z/w)'<0.
\end{equation}
Moreover,
\[
   z   (\frac{\pi}{2})=w'(\frac{\pi}{2}),\qquad   z   '(\frac{\pi}{2})=0,\qquad
 \lim_{r\to0^+}\mathcal W_{z}(r)=\lim_{r\to0^+}\Phi(r)\mathcal W_w(r)=0,
\]
and $\Phi=O((\frac{\pi}{2}-r)^{1-b})$ as $r\to \frac{\pi}{2}^-$.
\end{lemma}
\begin{proof}

\emph{The case $k\ge4$.}
Take
\[
 \begin{gathered}
   z   =\sqrt{\sin r}\left(w'+\frac1{b}\cot r\,w\right),\qquad
 \Phi=\frac{C}{d}\cos r\,\sin^{-5/2}r\,w^{-b},\\
 C=\frac{(2k-5)(2k-3)}{16}>0.
 \end{gathered}
\]
Using $\mathcal Dw=dw^p$ and
$\mathcal Dw'=(k-1)\csc^2r\,w'$, after a direct calculation one obtains
\[
 \mathcal D z=-C\frac{\cos r}{\sin^{5/2}r}w.
\]
The Wronskian identity now yields $\mathcal W_{z}'=\Phi\mathcal W_w'$. Furthermore,
\[
 \tan r\,\frac{\Phi'}\Phi=\xi-t-\frac52<0
\]
by Lemma~\ref{lem:phase-bound}. Since
\[
 \frac{z}{w}=\frac{\cos r}{b\sqrt{\sin r}}(1-\xi),
\]
equation~\eqref{eq:xi} gives
\[
 \left(\frac{z}{w}\right)'=
 -\frac{\cos r}{b\sqrt{\sin r}\tan r}
 \left(t+\frac{1+\xi}{2}+\frac{\xi G}{b}\right)<0.
\]
At zero, $\mathcal W_{z},\Phi\mathcal W_w=O(r^{k-5/2})\to0$.
The endpoint assertions at $\frac{\pi}{2}$ follow from \eqref{eq:mid-expansion}
and $w(\frac{\pi}{2})=w''(\frac{\pi}{2})=0$.

\emph{The cases $k=2,3$.}
Put $\gamma=(k-2)/2$ and take
\[
   z   =\sin r\,w'+\gamma\cos r\,w,\qquad
 \Phi=\cos r(\mathfrak a-\mathfrak bw^{-b}),\qquad
 \mathfrak a=\frac{k-1}{2},\quad \mathfrak b=\frac3{8d}.
\]
Again, a direct calculation gives
\[
 \mathcal D z=\cos r\left(\frac38w
                    -\frac{d(k-1)}2w^p\right),
\]
and hence $\mathcal W_{z}'=\Phi\mathcal W_w'$. Lemma~\ref{lem:convex} implies
$-\cot r\,w'/w>1$, so
\[
 \Phi'=\sin r\left[-\mathfrak a+\mathfrak bw^{-b}
           \left(1+b\cot r\,\frac{w'}w\right)\right]<0,
\]
because $b>1$. Finally,
\[
 \frac{z}{w}=\cos r\left(\gamma-\frac{\xi}{b}\right),\qquad
 \left(\frac{z}{w}\right)'=
 -\frac{\cos r}{\tan r}
 \left[\gamma t+\frac{\xi}{b}\left(1+\frac G{b}\right)\right]<0.
\]
Here we used $G>0$ from Lemma~\ref{lem:phase}. At zero,
$\mathcal W_{z}=O(r^k)$, $\mathcal W_w=O(r^k)$ and $\Phi=O(1)$.
The assertions at $\frac{\pi}{2}$ follow as in the preceding case.
\end{proof}

At this point we can prove an important lemma, which gives the zero pattern of $\varphi$ and the positivity of $\mathcal W_w$, similarly to \cite[Lemma 4.2]{Tang}.
\begin{lemma}\label{lem:shooting}
For every $k\ge2$, the solution $\varphi$ of \eqref{eq:shoot} has exactly
one zero $\tau$ in $[0,\frac{\pi}{2}]$, and $\tau\in(0,r_0)$, where $r_0\in(0,\frac{\pi}{2})$ is the unique point with $w(r_0)=1$. Moreover,
\[
 \varphi>0\ \text{on }[0,\tau),\qquad
 \varphi<0\ \text{on }(\tau,\frac{\pi}{2}],\qquad
 \mathcal W_w>0\ \text{on }(0,\frac{\pi}{2}].
\]
In particular,  $\varphi\left(\frac\pi2\right)<0.$
\end{lemma}
\begin{proof}
\emph{Step 1. The first zero occurs before $r_0$.}
On $(0,r_0)$,
\[
 \mathcal D(f(w))=f''(w)(w')^2>0.
\]
Thus, with
$P=a(f(w)\varphi'-f'(w)w'\varphi)$,
\[
 P'=-af''(w)(w')^2\varphi,\qquad P(0)=0.
\]
If $\varphi$ had no zero in $(0,r_0)$, then $P(r_0)<0$.
On the other hand, $\varphi(r_0)\ge0$ and $f'(1)=d>0$
would give
\[
 P(r_0)=-a(r_0)dw'(r_0)\varphi(r_0)\ge0,
\]
a contradiction. Hence the first zero $\tau$ lies in $(0,r_0)$,
it is simple, and $\varphi'(\tau)<0$.

\emph{Step 2. Positivity of $\mathcal W_w$ in the open interval.}
Use the functions from Lemma~\ref{lem:radial-comparison} and define
\[
 \mathcal F=\mathcal W_{z}-\Phi\mathcal W_w.
\]
By \eqref{eq:comparison-properties},
\begin{equation}\label{eq:Fprime}
 \mathcal F'=-\Phi'\mathcal W_w,\qquad \lim_{r\to0^+}\mathcal F(r)=0.
\end{equation}
Suppose that $\mathcal W_w$ has a first zero $s\in(0,\frac{\pi}{2})$. Since $\mathcal W_w>0$ near zero,
we have $\mathcal W_w>0$ on $(0,s)$, and therefore $\mathcal F(s)>0$.
As $\mathcal W_w(s)=0$, it follows that $\mathcal W_{z}(s)>0$.
But $\mathcal W_w'(s)\le0$ and \eqref{eq:Wwprime} imply $\varphi(s)\le0$.
Equality would force $\varphi'(s)=0$ as well, contrary to uniqueness
for the linear equation. Thus $\varphi(s)<0$ and
$\varphi'(s)/\varphi(s)=w'(s)/w(s)$. Consequently
\[
 \mathcal W_{z}(s)=-a(s)w(s)\varphi(s)
                    \left(\frac{z}{w}\right)'(s)<0,
\]
a contradiction. This proves $\mathcal W_w>0$ on $(0,\frac{\pi}{2})$ and, by
\eqref{eq:Fprime}, $\mathcal F'>0$ there.

\emph{Step 3. There is no second interior zero.}
At a second zero $\sigma\in(\tau,\frac{\pi}{2})$ we would have
$\varphi'(\sigma)>0$ and
$\mathcal W_w(\sigma)=-a(\sigma)w(\sigma)\varphi'(\sigma)<0$.
Hence $\varphi<0$ on $(\tau,\frac{\pi}{2})$ and $\varphi(\frac{\pi}{2})\le0$.

\emph{Step 4. The midpoint value is strictly negative.}
Suppose that $\varphi(\frac{\pi}{2})=0$. Writing $\varrho=\frac{\pi}{2}-r$, regularity gives
$w=O(\varrho)$, $\varphi=O(\varrho)$, and hence, by \eqref{eq:Wwprime},
$\mathcal W_w'=O(\varrho^{p+1})$. Since $\mathcal W_w(\frac{\pi}{2})=0$, we obtain
\[
 \mathcal W_w=O(\varrho^{p+2}),\qquad
 \Phi\mathcal W_w=O(\varrho^{p+3-b})=O(\varrho^4)\longrightarrow0.
\]
Lemma~\ref{lem:radial-comparison} therefore yields
\[
 \mathcal F(\frac{\pi}{2})=\mathcal W_{z}(\frac{\pi}{2})=w'(\frac{\pi}{2})\varphi'(\frac{\pi}{2}).
\]
The left side is strictly positive by \eqref{eq:Fprime} and $\mathcal W_w>0$.
Thus $\varphi'(\frac{\pi}{2})<0$, which is impossible if $\varphi<0$ immediately
before $\frac{\pi}{2}$ and $\varphi(\frac{\pi}{2})=0$. We conclude that $\varphi(\frac{\pi}{2})<0$.
Finally, $\mathcal W_w(\frac{\pi}{2})=w'(\frac{\pi}{2})\varphi(\frac{\pi}{2})>0$.
\end{proof}


\begin{lemma}\label{lem:midpoint}
For every $k\ge2$, the solution $\varphi$ of \eqref{eq:shoot} satisfies $\varphi'(\frac{\pi}{2})<0$.
\end{lemma}
\begin{proof}
Use $z,\Phi,\mathcal W_{z}$ from Lemma~\ref{lem:radial-comparison}.
For $k\ge4$, $\Phi\mathcal W_w$ vanishes at both endpoints. Moreover,
$\Phi'\mathcal W_w=O(r^{k-7/2})$ near zero and $\Phi'\mathcal W_w=O((\frac{\pi}{2}-r)^{-b})$ near $\frac{\pi}{2}$,
so the integral below converges. By Lemma~\ref{lem:shooting},
\[
 w'(\frac{\pi}{2})\varphi'(\frac{\pi}{2})=\mathcal W_{z}(\frac{\pi}{2})
 =-\int_0^{\frac{\pi}{2}} \Phi'\mathcal W_w\,dr>0.
\]
For $k=2,3$, the weight $\Phi$ tends to $-\infty$ at $\frac{\pi}{2}$, but
$\Phi\mathcal W_w'=O((\frac{\pi}{2}-r)^2)$ is integrable. Also $\mathfrak a\ge \mathfrak b>0$ and
$w(\tau)>1$, so $\Phi(\tau)>0$. Using the sign change of $\mathcal W_w'$
at $\tau$ from Lemma~\ref{lem:shooting}, we obtain
\[
 \mathcal W_{z}(\frac{\pi}{2})=\Phi(\tau)\mathcal W_w(\frac{\pi}{2})
       +\int_0^{\frac{\pi}{2}}[\Phi-\Phi(\tau)]\mathcal W_w'\,dr>0.
\]
The integrand is strictly positive on both sides of $\tau$.
Since $\mathcal W_{z}(\frac{\pi}{2})=w'(\frac{\pi}{2})\varphi'(\frac{\pi}{2})$ and $w'(\frac{\pi}{2})<0$,
the assertion follows in both cases.
\end{proof}

Now we can prove radial nondegeneracy, by which we mean nondegeneracy in the $\Gamma$-invariant subspace $H^1_\Gamma(\s^{2k-1})$,  with $\Gamma=O(k)\times O(k)$.

\begin{proposition}\label{prop:radial}
For every $k\ge2$, the only regular
solution of \eqref{eq:mode} with $j=\ell=0$ is zero.
\end{proposition}
\begin{proof}
Every regular solution at $0$ must be a multiple of
$\varphi$. Since $f'(w)$ is invariant under $r\mapsto\pi-r$, a
nonzero solution of \eqref{eq:moderad} must be either even or odd
about $\frac{\pi}{2}$. The even case requires $\varphi'(\frac{\pi}{2})=0$, contrary to
Lemma~\ref{lem:midpoint} and the odd case requires $\varphi(\frac{\pi}{2})=0$, contrary
to Lemma~\ref{lem:shooting}.
\end{proof}

Besides nondegeneracy, Lemma~\ref{lem:shooting} also yields the following uniqueness result (see also \cite[Theorem~1]{SW}).

\begin{lemma}\label{lem:unique-profile}
For every $k\ge2$, there exists exactly one solution $w$ of \eqref{eq:w} satisfying
\eqref{eq:solution}, equivalently, the symmetric single-node solution is unique.
\end{lemma}

\begin{proof}
Existence follows from the construction in \cite{FP}. The energy
$E$ defined in \eqref{energy} is nonincreasing on $(0,\frac{\pi}{2})$, so the solution $w(\cdot;\beta)$ of \eqref{ivp} is defined and bounded on $[0,\frac{\pi}{2}]$, and it depends
continuously on $\beta$, uniformly on $[0,\frac{\pi}{2}]$. Set
\begin{equation} \label{setcalS}
 g(\beta)=w\bigl(\tfrac{\pi}{2};\beta\bigr),\qquad
 \mathcal S=\bigl\{\beta>0:\ w(\cdot;\beta)>0\ \text{on }[0,\tfrac{\pi}{2}),
 \ g(\beta)=0\bigr\}.
\end{equation}
By the argument above (odd reflection and the energy
identity), $\beta\in\mathcal S$ if and only if $w(\cdot;\beta)$, extended
oddly across $\frac{\pi}{2}$, satisfies \eqref{eq:solution}.
Note that $\beta\in\mathcal S$ forces $\beta>1$, indeed, if $\beta\le1$ then
$E(0)<0$, whereas $E(\frac{\pi}{2})=\frac12 w'(\frac{\pi}{2};\beta)^2\ge0$, against the monotonicity of $E$.

For $\beta\in\mathcal S$, the derivative
$\partial_\beta w(\cdot;\beta)$ is the solution
$\varphi$ of \eqref{eq:shoot}, so Lemma~\ref{lem:shooting} gives
\begin{equation}\label{eq:g-transversal}
 g'(\beta)=\varphi\bigl(\tfrac{\pi}{2}\bigr)<0
 \qquad\text{for every }\beta\in\mathcal S.
\end{equation}

Every zero of $w(\cdot;\beta)$ in $[0,\frac{\pi}{2}]$ is simple and lies
in $(0,\frac{\pi}{2}]$, because $w(0)=\beta>0$. Let $N(\beta)$ be the
number of zeros in $(0,\frac{\pi}{2})$.

Suppose that $\beta_1<\beta_2$ both lie in $\mathcal S$.
By \eqref{eq:g-transversal}, $g<0$ on some interval $(\beta_1,\beta_1+\delta)$.
Since $w(0;\beta)>0$, this gives $N\ge1$ there.
By \eqref{eq:g-transversal} again, and by the positivity of
$w(\cdot;\beta_2)$ on $[0,\frac{\pi}{2})$ together with
$w'(\frac{\pi}{2};\beta_2)\ne0$, we have $N=0$ and $g>0$ on some interval
$(\beta_2-\delta,\beta_2)$.

Let 
\[
 {\hat\beta}=\sup\{\beta\in(\beta_1,\beta_2): N(\beta)\ge1\}
 \in[\beta_1+\delta,\beta_2-\delta].
\]
Then $N=0$ on $({\hat\beta},\beta_2)$, so $g\ge0$ there. By continuity,
$w(\cdot;{\hat\beta})\ge0$ on $[0,\frac{\pi}{2}]$. A nonnegative solution
with an interior zero would have a double zero, so
$w(\cdot;{\hat\beta})>0$ on $(0,\frac{\pi}{2})$.

On the other hand, there are $\beta_j\uparrow{\hat\beta}$ and zeros
$r_j\in(0,\frac{\pi}{2})$ of $w(\cdot;\beta_j)$. A limit point $\bar r$
of the $r_j$ satisfies $w(\bar r;{\hat\beta})=0$. Uniform convergence and
$w(0;{\hat\beta})>0$ exclude $\bar r=0$, so $\bar r=\frac{\pi}{2}$.

Thus ${\hat\beta}\in\mathcal S$, and \eqref{eq:g-transversal} gives $g<0$
immediately to the right of ${\hat\beta}$. This contradicts $g\ge0$ on
$({\hat\beta},\beta_2)$.
\end{proof}
\section{The modes \texorpdfstring{$(1,0)$ and $(0,1)$}{(1,0) and (0,1)}}
\begin{lemma}\label{lem:degree-one}
The regular solutions of \eqref{eq:mode} in degrees $(1,0)$ and $(0,1)$
are spanned, respectively, by $h_{10}$ and $h_{01}$ defined in
\eqref{eq:geometric-modes}.
\end{lemma}
\begin{proof}
Substitution in \eqref{eq:mode} verifies the formulas. In degree $(1,0)$, the indicial roots of \eqref{eq:mode} are $0$ and $-(k-2)$ at $r=0$, and $1$ and $-(k-1)$ at $r=\pi$; for $k=2$ the root $0$ at $r=0$ is double and the second solution behaves like $\log r$. Hence the solutions that are bounded near $r=0$ form a one-dimensional space, which contains $h_{10}$. The same argument applies to $h_{01}$ in degree $(0,1)$, with the roles of the two endpoints exchanged. Thus these nonzero regular
solutions span the corresponding kernels. Notice that
$h_{10}(\pi-r)=-h_{01}(r)$.
\end{proof}

\section{The mode \texorpdfstring{$(1,1)$}{(1,1)}}
\begin{lemma}\label{lem:ground}
The regular kernel in the mode $(1,1)$ is spanned by $w'$. For regular $h$ vanishing at both endpoints,
\begin{equation}\label{eq:ground}
 \int_0^\pi a\left[h'^2-f'(w)h^2
                    +\frac{k-1}{\sin^2r}h^2\right]dr\ge0.
\end{equation}
\end{lemma}
\begin{proof}
Set $\Theta=-w'$. By \eqref{eq:solution},
$\Theta>0$ on $(0,\pi)$. Differentiating \eqref{eq:w} shows that $w'$ solves \eqref{eq:mode}
for $j=\ell=1$, since $\lambda_1=k-1$, therefore
\begin{equation}\label{eq:theta}
 (a\Theta')'
 =a\left(\frac{k-1}{\sin^2 r}-f'(w)\right)\Theta.
\end{equation}
For any admissible $h$, we have the following identity
\begin{align*}
 a\Theta^2
       \left(\left(\frac{h}{\Theta}\right)'\right)^2 = a\left(h'-\frac{\Theta'}{\Theta}h\right)^2
 &=a(h')^2-a\Theta'\left(\frac{h^2}{\Theta}\right)'\\
 &=a\left[(h')^2-f'(w)h^2
       +\frac{k-1}{\sin^2 r}h^2\right]
   -\left(a\frac{\Theta'}{\Theta}h^2\right)',
\end{align*}
where the last equality follows from \eqref{eq:theta}.
Since $\Theta$ vanishes linearly at $0$ and $\pi$ and a regular
coefficient in the $(1,1)$ mode satisfies
$h(r)=O(r)$ as $r\to0$ and $h(r)=O(\pi-r)$ as $r\to\pi$,
we have
\[
 \left.a\frac{\Theta'}{\Theta}h^2\right|_{r\to0}=0,
 \qquad
 \left.a\frac{\Theta'}{\Theta}h^2\right|_{r\to\pi}=0.
\]
Integration therefore yields
\begin{equation}\label{eq:picone}
 \int_0^\pi a\left[(h')^2-f'(w)h^2
           +\frac{k-1}{\sin^2 r}h^2\right]dr
 =\int_0^\pi a\Theta^2
       \left(\left(\frac{h}{\Theta}\right)'\right)^2dr
 \ge 0.
\end{equation}
If $h$ is a regular solution to \eqref{eq:mode} with $j=\ell=1$, multiplication of that equation by $ah$, followed by integration by parts, gives an equality in \eqref{eq:picone}. Hence $h$ is a multiple of
$\Theta=-w'$. 

\end{proof}

\section{The modes \texorpdfstring{$(j,\ell)\ne(1,1)$ with $j,\ell\ge1$}{(j,l) other than (1,1) with j,l>=1}}
\begin{lemma}\label{lem:mixed}
If $j,\ell\ge1$ and $(j,\ell)\ne(1,1)$, the regular kernel of
\eqref{eq:mode} is trivial for every $k\ge2$.
\end{lemma}
\begin{proof}
Multiplication of \eqref{eq:mode} by $ah$ and integration give
\[
 0=\int_0^\pi a\left[h'^2-f'(w)h^2
                  +\frac{k-1}{\sin^2r}h^2\right]dr
   +\int_0^\pi a\left[
       \frac{\lambda_j-\lambda_1}{4\cos^2(r/2)}
      +\frac{\lambda_\ell-\lambda_1}{4\sin^2(r/2)}\right]h^2\,dr.
\]
The first integral is nonnegative by \eqref{eq:ground}. The second is
strictly positive unless $h=0$, since at least one degree exceeds one.
\end{proof}

\section{The modes \texorpdfstring{$(0,\ell)$ or $(\ell,0)$ with $\ell\ge2$ for large $k$}{(0,l) or (l,0) with l>=2 for large k}}\label{sec:one-sided}

We consider the eigenvalue problem
\begin{equation}\label{eq:onesided}
 \LL h+\Lambda\rho h=0,\qquad
 \rho(r)=\frac{a(r)}{4\sin^2(r/2)}.
\end{equation}
In the usual Sturm--Liouville convention the eigenvalue parameter is
$-\Lambda$. Thus the order of the spectral values is reversed when they
are expressed in terms of $\Lambda$, so the first Sturm eigenfunction
corresponds to the largest value of $\Lambda$.

The eigenfunction $\zeta:=h_{01}$ solves \eqref{eq:onesided} at
$\Lambda=\lambda_1=k-1$. We will show that $\zeta$ has exactly one zero
in $(0,\pi)$. Sturm oscillation will then identify it as the second
eigenfunction. In particular, there is exactly one spectral value
$\Lambda_*>\lambda_1$, and its eigenfunction can be chosen strictly
positive. Therefore, if a mode of degree $\ell\ge2$ exists, its
parameter $\lambda_\ell$ must equal this single value $\Lambda_*$.

It is convenient to compare $\Lambda_*$ with the discrete sequence
$\{\lambda_\ell=\ell(\ell+k-2),\  \ell=0,1,2,\ldots \}$ by introducing a continuous parameter $\alpha>0$ through
\[
 \Lambda=\alpha(\alpha+k-2).
\]
The map $\alpha\mapsto\alpha(\alpha+k-2)$ is strictly increasing for
$\alpha>0$, and for integer $\alpha=\ell$ it gives exactly
$\Lambda=\lambda_\ell$. Thus $\alpha$ is only a reparametrization of
$\Lambda$ and it allows us to ask whether the principal value $\Lambda_*$ corresponds to an integer $\alpha$. Since $\Lambda_*>\lambda_1$, there is a unique
$\alpha_*>1$ such that
\[
 \Lambda_*=\alpha_*(\alpha_*+k-2).
\]
Consequently, a nonzero mode $(0,\ell)$ with $\ell\ge2$ can occur only
if $\ell=\alpha_*$. The rest of this section excludes this possibility
for sufficiently large $k$.

\begin{lemma}\label{lem:principal}
The function $\zeta:=h_{01}$ has exactly one zero in $(0,\pi)$.
There is exactly one spectral value $\Lambda_*>\lambda_1$ in
\eqref{eq:onesided}. Its eigenfunction $h_*$ can be chosen positive.
Writing $\Lambda_*=\alpha_*(\alpha_*+k-2)$, one has $\alpha_*>1$.
Consequently, a nonzero mode $(0,\ell)$ with $\ell\ge2$ must satisfy
$\ell=\alpha_*$.
\end{lemma}
\begin{proof}
On $(0,\frac{\pi}{2}]$, the explicit formula
\[
 \zeta(r)=\frac{2k-3}{2}w(r)\sin\frac r2
                 -2w'(r)\cos\frac r2
\]
and the assumptions $w>0$, $w'<0$ imply $\zeta>0$. It remains
to count its zeros on $(\frac{\pi}{2},\pi)$.

For $r$ in this interval, set $s=-\cos r\in(0,1)$ and recall
$W(x)=w(\arccos x)$. The antisymmetry of $w$ implies that $W$ is odd
and $W'$ is even. Hence
\[
 w(r)=W(-s)=-W(s),\qquad
 w'(r)=-\sin r\,W'(-s)=-\sin r\,W'(s).
\]
Using also $\cos^2(r/2)=(1-s)/2$ and
$\sin r=2\sin(r/2)\cos(r/2)$, substitution into the formula for
$\zeta$ gives
\begin{equation}
 \zeta(r)=\sin\frac r2\,D(s),\qquad
 D(s):=2(1-s)W'(s)-\frac{2k-3}{2}W(s).
\end{equation}
Because $\sin(r/2)>0$, the zeros of $\zeta$ on this half-interval are
exactly the zeros of $D$ in $(0,1)$. At the endpoints,
\[
 D(0)=2W'(0)>0,\qquad
 D(1)=-\frac{2k-3}{2}W(1)<0.
\]
Thus $D$ has at least one zero.

To show that every zero is crossed downwards, write
\[
 H(s):=b(1-s)\frac{W'(s)}{W(s)},\qquad
 D(s)=\frac{2W(s)}{b}\bigl(H(s)-1\bigr),
\]
where $b=4/(2k-3)$. For $r=\arccos s$ in $(0,\frac{\pi}{2})$, $H$ satisfies
\[
 H(s)=\frac{s}{1+s}\,\xi(r).
\]
At a zero $s_0$ of $D$, we therefore have $H(s_0)=1$. Differentiating
this identity and using \eqref{eq:xi}, together with
$t=\tan^2r=(1-s^2)/s^2$, yields
\[
 H'(s_0)=-\frac{1+s_0G(\arccos s_0)/b}{1-s_0^2}<0,
\]
where the strict inequality follows from $G>0$ in
Lemma~\ref{lem:phase}. Since
\[
 D'(s_0)=\frac{2W(s_0)}{b}H'(s_0),
\]
every zero of $D$ is crossed strictly downwards. Together with the
endpoint signs, this proves that $D$, and hence $\zeta$, has exactly
one zero on the second half of the interval. Thus $\zeta$ has exactly
one zero in $(0,\pi)$.

We now justify the Sturm ordering in the presence of the singular endpoint
$r=0$. Since $\LL h=-(ah')'-af'(w)h$, equation
\eqref{eq:onesided} is the generalized eigenvalue problem
\[
 -(ah')'-af'(w)h=\nu\rho h,\qquad \nu:=-\Lambda.
\]
Introduce the change of coordinates and the positive factor
\[
 \varsigma=\log\tan\frac r4,
 \qquad
 \mathfrak g(r)=(a(r)\rho(r))^{1/4},
 \qquad
 \chi(\varsigma)=\mathfrak g(r)h(r).
\]
Then $d\varsigma/dr=\sqrt{\rho/a}$, and the equation becomes
\[
 -\chi_{\varsigma\varsigma}
 +\mathcal V(\varsigma )\chi=\nu \chi,
 \qquad
 \mathcal V (\varsigma )
 =-4\sin^2\frac r2\,f'(w(r))
   +\frac{\mathfrak g_{\varsigma \varsigma }}{\mathfrak g}.
\]
Here the derivatives of $\mathfrak g$ are taken with respect to $\varsigma $. As $r\to0$, one has
$\varsigma \to-\infty$ and
$\mathfrak g(r)\sim C r^{(k-2)/2}$ for some $C>0$. Since $f'(w(r))$ stays bounded,
\[
 \lim_{\varsigma \to-\infty}\mathcal V (\varsigma )
 =\frac{(k-2)^2}{4}\ge0.
\]
At $r=\pi$ we impose the regular condition, namely that $h$ is bounded and $h'(r)=O(\pi-r)$. In these coordinates we get
\[
 \chi(\varsigma)=C(-\varsigma)^{(k-1)/2}(1+O(\varsigma^2))\qquad(\varsigma\to0^-),
\]
The regularity condition at $r=\pi$ selects the Friedrichs
realization of the transformed operator at $\varsigma=0$.
In particular, when $k=2$, it excludes the logarithmic branch
$(-\varsigma)^{1/2}\log(-\varsigma)$ associated with the
critical inverse-square potential.

At $\varsigma=-\infty$, we impose square integrability.
Thus the essential spectrum starts at a nonnegative value, and the
negative eigenvalues $\nu=-\Lambda$ under consideration are isolated.
The standard Sturm oscillation ordering therefore applies (equivalently,
one may obtain it by regular endpoint exhaustion), see, e.g.,~\cite{Weidmann,Zettl}.

Since $\zeta$ has exactly one interior zero, it is the second
eigenfunction. Its Sturm eigenvalue is $-\lambda_1$, so exactly one
eigenvalue lies below it, the corresponding first eigenfunction is
strictly positive and has parameter $\Lambda_*>\lambda_1$. All remaining
spectral values have parameter less than $\lambda_1$. The strict
monotonicity of $\alpha\mapsto\alpha(\alpha+k-2)$ now gives the stated
$\alpha_*>1$ and the conclusion for every possible mode $(0,\ell)$ with
$\ell\ge2$.
\end{proof}

Factoring the regular power as
$h(r)=\sin^\alpha(r/2)y(\cos r)$ gives
\begin{equation}\label{eq:alpha}
 (1-x^2)y''-[\alpha+(k+\alpha)x]y'
 +\left(V-\frac{\alpha(\alpha+2k-2)}4\right)y=0.
\end{equation}

\begin{lemma}\label{lem:alpha2}
For every $k\ge2$, $\alpha_*\ne2$.
\end{lemma}
\begin{proof}
If $\alpha_*=2$, the positive principal eigenfunction gives a positive
regular $y$ on $[-1,1]$ satisfying
\[
 (1-x^2)y''-[2+(k+2)x]y'+(V-k)y=0.
\]
For $0\le x\le1$, write $\bar y(x)=y(-x)$ and set
\[
 B=(1-x^2)(y'\bar y-y\bar y')-2y\bar y.
\]
The equations of $y$ and $\bar y$, using the evenness of $V$, give
\[
 B'=\frac{kx}{1-x^2}(B+2y\bar y).
\]
Since $B(1)=-2y(1)y(-1)<0$ and $B'>0$ at every interior zero,
$B<0$ on $(0,1)$. Integrating the resulting logarithmic inequality gives
\begin{equation}\label{eq:reflection}
 \frac{y(x)}{y(-x)}<\frac{1+x}{1-x},\qquad0<x<1.
\end{equation}
The integrating factor for the equation of $y$ is
\[
 p_2(x)=(1-x)^{(k+4)/2}(1+x)^{k/2},\qquad
 \frac{p_2'}{p_2}=-\frac{2+(k+2)x}{1-x^2}.
\]
Equation~\eqref{eq:Y} shows that substitution of $Y$ in the equation
of $y$ gives the residual $-2Y'$. Its Wronskian identity with $y$,
with vanishing endpoint terms, therefore implies
\[
 0=\int_{-1}^1\frac{p_2yY'}{1-x^2}\,dx
  =\int_0^1\frac{Y'(x)}{1-x^2}
       [p_2(x)y(x)-p_2(-x)y(-x)]\,dx.
\]
But \eqref{eq:reflection} and
$p_2(x)/p_2(-x)=((1-x)/(1+x))^2$ make the bracket strictly negative.
Since $Y'>0$ on $(0,1)$ by Lemma~\ref{lem:convex}, the last integral
is negative, a contradiction.
\end{proof}

\begin{lemma}\label{lem:amplitude}
There exists $k_1$ such that, for every $k\ge k_1$,
\begin{equation}\label{eq:amplitude}
 w(0)^{b}<2.
\end{equation}
\end{lemma}
\begin{proof}
Suppose, for a contradiction, that this fails, and consider an unbounded sequence of integers $k$
for which $S(0)\ge2$. All constants below are independent of this sequence.
We use the auxiliary functions in \eqref{eq:phase-def}.

\emph{Step 1. The first hitting point of $\xi=1$ stays away from $\pi/2$.}
We first record bounds that hold as long as $\xi<1$. By the expansion at
the origin, $\mu$ extends continuously to $r=0$ with $\mu(0)=k>d$.
If $\mu$ reached $d$ for the first time at some $r$ while $\xi<1$, then
$\mu'(r)\le0$. But, evaluating \eqref{eq:mu} at $\mu=d$ and using
$k-d=p/b$ and $d=bc$, we get
\[
 \tan r\,\mu'=d\left(k-d-\frac p{b}\xi\right)
              =pc(1-\xi)>0,
\]
a contradiction. Thus $\mu>d$ while $\xi<1$. Lemma~\ref{lem:phase} also
gives $\mu<k$. Since $\mu\xi=dt(S-1)$, we have $S>1$ and
\begin{equation}\label{eq:sandwich}
 \frac{d}{k}\,t(S-1)<\xi<t(S-1),\qquad t=\tan^2r,
\end{equation}
throughout this interval.

We now use $t$ as the independent variable. Differentiating
$S=w^b$ and using the definition of $\xi$ gives
\[
 \frac{d}{dt}\log S=-\frac{\xi}{2t(1+t)}.
\]
The upper bound in \eqref{eq:sandwich} then implies
\[
 \frac{d}{dt}\log\frac{S-1}{S}
 =-\frac{\xi}{2t(1+t)(S-1)}
 >-\frac{1}{2(1+t)}.
\]
Integrating from $0$ to $t$ and using $S(0)\ge2$, hence
$(S(0)-1)/S(0)\ge1/2$, yields
\[
 \frac{S(t)-1}{S(t)}>
 \frac{1}{2\sqrt{1+t}},
 \qquad
 S(t)-1>\frac{1}{2\sqrt{1+t}-1}.
\]
Combining this with the lower bound in \eqref{eq:sandwich}, we obtain
\[
 \xi(t)>\frac{d}{k}\,
          \frac{t}{2\sqrt{1+t}-1}
 =\frac{2k-1}{4k}
          \frac{t}{2\sqrt{1+t}-1}.
\]
Choose a fixed $t_*>0$ such that
\[
 \frac{t_*}{2\sqrt{1+t_*}-1}>2,
 \qquad
 \eta_1:=\arctan\sqrt{t_*}<\frac{\pi}{2}.
\]
If $\xi$ had not reached $1$ by $r=\eta_1$, the preceding estimate at
$t=t_*$ would give
\[
 \xi(\eta_1)>\frac{2k-1}{4k}
                  \frac{t_*}{2\sqrt{1+t_*}-1}>1
\]
for all sufficiently large $k$, a contradiction. Since $\xi$ extends
continuously to $\xi(0)=0$, it follows that its first hitting point
$r_{1,k}$ exists and satisfies $r_{1,k}<\eta_1$ for all sufficiently
large $k$.

\emph{Step 2. The first hitting point of $\xi=2$ also stays in a fixed interval.}
At $r_{1,k}$, we have $\xi=1$. Since $G>0$,
\[
 \mu(r_{1,k})>k-1-\frac1b=d.
\]
The identity $\mu\xi=dt(S-1)$ therefore gives
\[
 t(r_{1,k})\bigl(S(r_{1,k})-1\bigr)>1,
 \qquad\text{equivalently}\qquad
 S(r_{1,k})\sin^2r_{1,k}>1.
\]
Moreover, differentiating $S\sin^2r$ gives
\[
 \frac{d}{dr}\log(S\sin^2r)=(2-\xi)\cot r.
\]
Thus $S(r)\sin^2r$ is strictly increasing as long as $\xi<2$, so the
same strict inequality holds up to the first point where $\xi$ reaches
$2$.

It remains to show that this point occurs before a bound independent of
$k$. By \eqref{eq:xi} and $G>0$,
\[
 \frac{\xi'}{\xi}
 =\cot r\left(1+t+\frac{G}{b}\right)
 >(1+t)\cot r
 =\frac{d}{dr}\log\tan r.
\]
Integrating from $r_{1,k}$, where $\xi(r_{1,k})=1$, shows that while
$\xi<2$,
\[
 \xi(r)>\frac{\tan r}{\tan r_{1,k}}.
\]
Consequently $\xi$ must reach $2$ before
$\tan r=2\tan r_{1,k}$. If $r_2$ denotes its first hitting point, then
\[
 r_{1,k}<r_2<\eta_2:=\arctan(2\tan\eta_1)<\frac{\pi}{2},
 \qquad S(r_2)\sin^2r_2>1.
\]
At $r_2$, the identity \eqref{eq:muxi} gives
$2\mu(r_2)=dt(r_2)(S(r_2)-1)>d$, hence $\mu(r_2)>d/2$. Using
$G=b\mu+\xi-b(k-1)$ and $b(k-1)=2+b/2$, we conclude that
\begin{equation}\label{eq:gap2}
 G(r_2)>\frac{bd}{2}+2-b(k-1)
           =\frac12-\frac{b}{4}.
\end{equation}
In particular, for every fixed $\varepsilon_0\in(0,1/2)$,
$G(r_2)>\varepsilon_0$ once $k$ is sufficiently large.

\emph{Step 3. Persistence of positivity of $G$.}
Fix $\eta_3\in(\eta_2,\frac{\pi}{2})$ and $\varepsilon_0\in(0,1/2)$.
We claim that $G\ge\varepsilon_0$ on $[r_2,\eta_3]$ for large $k$.
The initial inequality follows from \eqref{eq:gap2}. Equations
\eqref{eq:xi}--\eqref{eq:mu} give the identity
\begin{equation}\label{eq:G-full}
 \begin{aligned}
 \tan r\,G'={}&2+\frac{b}{2}+G+(1+G)t
       -\left(2+\frac{b}{2}+G-\xi\right)\xi\\
 &+\frac G{b}\left(2\xi-2-\frac{b}{2}-G\right).
 \end{aligned}
\end{equation}
At a first downward contact with $G=\varepsilon_0$, monotonicity of
$\xi$ and \eqref{eq:phase-bound} imply
$2\le\xi\le\tan^2\eta_3+5/2$. The first line of \eqref{eq:G-full}
is bounded below independently of $k$, whereas the second is at least
$\varepsilon_0/b$ for large $k$. Thus $G'>0$, a contradiction.

\emph{Final step.}
From \eqref{eq:xi} we obtain, on $[r_2,\eta_3]$,
\[
 \frac{\xi'}\xi\ge\frac{\varepsilon_0}{b}\cot r.
\]
Consequently
\[
 \xi(\eta_3)\ge
 2\exp\left(\frac{\varepsilon_0}{b}
                 \log\frac{\sin\eta_3}{\sin\eta_2}\right)\longrightarrow\infty,
\]
contradicting $\xi(\eta_3)<\tan^2\eta_3+5/2$ from
Lemma~\ref{lem:phase-bound}. This proves \eqref{eq:amplitude}.
\end{proof}

Set $q=Y'/Y$ on $[0,1]$. Lemma~\ref{lem:convex} gives
$q(0)=0$ and $q>0$ on $(0,1)$. By \eqref{eq:Y},
\begin{equation}\label{eq:riccati}
 (1-x^2)(q'+q^2)-(k+2)xq+V-k=0.
\end{equation}

\begin{lemma}\label{lem:q}
For all sufficiently large $k$,
\begin{equation}\label{eq:qbound}
 \max_{[0,1]}q<\frac k6.
\end{equation}
\end{lemma}
\begin{proof}
Suppose that $Q=\max q\ge k/6$ along an unbounded sequence of integers $k$.
Since $q(0)=0$, a maximum point lies in $(0,1]$.

\emph{Interior maximum.}
Let $x_*\in(0,1)$ be a maximum point. We first show that
$\mu(\arccos x_*)>0$ and $x_*>\eta_0:=\sqrt{10}-3$.
Equation \eqref{eq:W} and the definition of $\mu$ in \eqref{eq:phase-def} give, with $r=\arccos x$,
\begin{equation}\label{eq:qmu}
 q(x)=\frac{x(k-\mu(r))}{1-x^2}.
\end{equation}
At $x=x_*$, differentiation of \eqref{eq:qmu} yields
\[
 \tan r\,\mu'(r)=-(k-\mu)(t+2).
\]
Combining this with \eqref{eq:mu} gives
\[
 0=\mu\left(k-2-\mu-\frac p{b}\xi\right)
       +(k-d)t+2k.
\]
If $\mu\le0$, positivity of $G$ implies
\[
 k-2-\mu-\frac p{b}\xi
 <k-2-\mu-p(k-1-\mu)=-1-b(k-1-\mu)<0.
\]
Every term in the preceding identity is then nonnegative and its last
term is positive, a contradiction. Thus $\mu>0$, and \eqref{eq:qmu}
gives $Q<kx_*/(1-x_*^2)$. Since $Q\ge k/6$, one obtains
$x_*^2+6x_*-1>0$, hence $x_*>\eta_0$.

Now, differentiation of \eqref{eq:riccati}, using
$q'=0$ and $q''\le0$, gives
\[
 V'(x_*)\ge(k+2)Q+2x_*Q^2.
\]
Also
\[
 V'=b(V+c)\frac{W'}W,\qquad
 W(x_*)=\int_0^{x_*}Y(s)\,ds
 \ge Y(x_*)\frac{1-e^{-Qx_*}}{Q}.
\]
The second inequality follows by integrating $Y'/Y\le Q$ backwards
from $x_*$. Combining these estimates yields
\begin{equation}\label{eq:maxineq}
 (k+2+2x_*Q)(1-e^{-Qx_*})\le b(V(x_*)+c).
\end{equation}
Lemma~\ref{lem:amplitude} gives the upper bound
\[
 b(V+c)=bcpW^{b}<2dp=k+O(1).
\]
Using the strict inequality $x_*>\eta_0$ proved above, the left side of
\eqref{eq:maxineq} is at least
\[
 (k+2+\eta_0k/3)(1-e^{-\eta_0k/6})
 =\left(1+\frac{\eta_0}{3}+o(1)\right)k,
\]
a contradiction.

\emph{Boundary maximum.}
Suppose that the maximum is at $1$. Evaluating \eqref{eq:Y} and its
derivative at $1$ gives
\begin{equation}\label{eq:endpoint-q}
 V(1)=k+(k+2)q(1),\qquad
 q'(1)=\frac{V'(1)-(k+2)q(1)-2q(1)^2}{k+4}.
\end{equation}
The equation \eqref{eq:W} at $1$ also yields
\[
 \frac{W'(1)}{W(1)}=\frac{V(1)-d}{kp},\qquad
 V'(1)=\frac4{k(2k+1)}(V(1)+c)(V(1)-d).
\]
By Lemma~\ref{lem:amplitude},
\[
 V(1)<c(2p-1)=\frac{(2k-1)(2k+5)}{16},
\]
so \eqref{eq:endpoint-q} implies $q(1)<k/4$.
Write $q(1)=\theta k$, where $1/6\le\theta<1/4$. Uniformly in
this range,
\[
 V'(1)=\left(2\theta^2+\frac{\theta}{2}\right)k^2+O(k),
\]
and consequently
\[
 V'(1)-(k+2)q(1)-2q(1)^2
 =-\frac{\theta}{2}k^2+O(k)<0
\]
for large $k$. Equation~\eqref{eq:endpoint-q} then gives $q'(1)<0$,
contrary to maximality at the right endpoint.
\end{proof}

\begin{lemma}\label{lem:alpha3}
For all sufficiently large $k$, $\alpha_*<3$.
\end{lemma}
\begin{proof}
Substitute $Y$ into \eqref{eq:alpha} with $\alpha=3$. By \eqref{eq:Y},
the residual is
\[
 R_3(x)=-(3+x)Y'(x)-\frac{2k+3}{4}Y(x).
\]
It is negative on $[0,1)$ by Lemma~\ref{lem:convex}. For $x=-s<0$,
evenness of $Y$ and Lemma~\ref{lem:q} give
\[
 \frac{R_3(-s)}{Y(s)}=(3-s)q(s)-\frac{2k+3}{4}
 <\frac k2-\frac{2k+3}{4}<0.
\]
Set $\Psi(r)=\sin^3(r/2)Y(\cos r)>0$. Returning to $r$,
\[
 \LL\Psi+3(k+1)\rho\Psi
 =-a\sin^3(r/2)R_3(\cos r)>0.
\]
If $\alpha_*\ge3$, then $\Lambda_*\ge3(k+1)$. With the positive
principal eigenfunction $h_*$ from Lemma~\ref{lem:principal}, define
$\mathcal Q=a(\Psi h_*'-h_*\Psi')$. Its derivative is
\[
 \mathcal Q'=h_*[\LL\Psi+3(k+1)\rho\Psi]
       +(\Lambda_*-3(k+1))\rho\Psi h_*>0.
\]
At $0$ we have $\Psi=O(r^3)$ and $h_*=O(r^{\alpha_*})$, at $\pi$ both
functions have vanishing derivatives. Thus
$\mathcal Q(0)=\mathcal Q(\pi)=0$, a contradiction.
\end{proof}

\begin{proposition}\label{prop:onesided}
For all sufficiently large $k$, the regular kernel of \eqref{eq:mode}
is trivial in degrees $(0,\ell)$ and $(\ell,0)$ for every $\ell\ge2$.
\end{proposition}
\begin{proof}
By Lemmas~\ref{lem:principal}, \ref{lem:alpha2}, and \ref{lem:alpha3},
a nontrivial solution in case $(0,\ell)$ would require an integer
$\ell=\alpha_*\in(1,3)\setminus\{2\}$, which is impossible.
By symmetry, the same holds for $(\ell,0)$.
\end{proof}

\begin{remark}\label{rem:converse}
Except for Lemma~\ref{lem:alpha3}, all the results of Sections~\ref{sec:separation}--\ref{sec:one-sided} hold for every $k\ge2$. By Lemma~\ref{lem:principal}, a nonzero mode $(0,\ell)$, or $(\ell,0)$, with $\ell\ge2$ can occur only if $\ell=\alpha_*$. Hence, if $\alpha_*(k)\notin\N$, the proof of Theorem~\ref{thm:main} applies without changes, and \eqref{eq:kernel} holds. 
\end{remark}

\begin{proof}[Proof of Theorem~\ref{thm:main}]
Proposition~\ref{prop:radial} excludes $(0,0)$. Lemmas~\ref{lem:ground} and
\ref{lem:degree-one} identify $(1,1)$, $(1,0)$ and $(0,1)$.
Lemma~\ref{lem:mixed} and Proposition~\ref{prop:onesided} exclude all
other possibilities for large $k$. Since $\dim\mathcal H_1=k$, it follows that  $\dim\ker L_v=k(k+2)$.
\end{proof}

\section{Geometric nondegeneracy in Euclidean space}\label{sec:euclidean-kernel}
Let $P_N=e_{2k}$ and let $\Pi:\s^n\setminus\{P_N\}\to\R^n$ be
stereographic projection
\[
 \Pi(y)=\frac{(y_1,\ldots,y_n)}{1-y_{n+1}},\qquad
 \Pi^{-1}(z)=\frac{(2z,|z|^2-1)}{1+|z|^2}.
\]
Here we use $y$ for coordinates in $\s^n$ and $z$ for coordinates in $\R^n$. Set
\begin{equation}\label{eq:stereographic}
 \sigma(z)=\left(\frac2{1+|z|^2}\right)^{(n-2)/2},\qquad
 u(z)=\sigma(z)v(\Pi^{-1}(z)).
\end{equation}
Conformal covariance shows that $u$ solves \eqref{eqRn} if and only if $v$ solves \eqref{eqSn}, see for example \cite[Lemma 1]{Ding}.
The same transformation identifies $\ker L_v$ in $H^1(\s^n)$ with the kernel of
\[
 L_u \phi=-\Delta \phi-\frac{n(n+2)}4|u|^{p-1}\phi
\]
in $D^{1,2}(\R^n)$ through the map
\begin{equation}\label{eq:transfer}
 \mathcal T[\phi](y)=\frac{\phi(\Pi(y))}{\sigma(\Pi(y))}.
\end{equation}
In particular, $\ker_{H^1(\s^n)} L_v$ and $\ker_{D^{1,2}(\R^n)} L_u$ have the same dimension.

Let $u$, given by \eqref{eq:stereographic}, be the \emph{symmetric single-node solution} of \eqref{eqRn}.
The choice of the pole $P_N=e_{2k}$ induces the decomposition
\[
 z=(z',z'')\in\R^k\times\R^{k-1}.
\]
We use the functions $Z_\alpha$ defined in
\eqref{eq:euclidean-basis}, whose indices run from $0$ to
$k(k+2)-1$ without repetition.
The first function is the dilation direction, the next $n$ functions
are translation directions, and the next $k$ are special conformal
directions, obtained by conjugating translations with the Kelvin
transform. The last $k(k-1)$ functions correspond to rotations in the
$(z_i,z_{k+a})$-planes. Each of these functions belongs to
$D^{1,2}(\R^n)$ and satisfies $L_uZ_\alpha=0$ by conformal covariance.

The selection of directions in \eqref{eq:euclidean-basis} reflects the
symmetries of $u$. Rotations within either Euclidean coordinate block
act trivially on $u$. Moreover, the remaining special conformal
directions are redundant and invariance under spherical rotations in
the $(y_{k+a},y_{2k})$-planes gives
\[
 -2z_{k+a}Z_0+|z|^2Z_{k+a}=Z_{k+a},
 \qquad 1\le a\le k-1.
\]
Consequently, the span of the functions in
\eqref{eq:euclidean-basis} coincides with the space
$\widetilde{\mathcal Z}_u$ generated by all infinitesimal conformal
invariances. Equality of this space with the full kernel is the
geometric nondegeneracy condition of Duyckaerts--Kenig--Merle~\cite{DKM}.

To verify linear independence and identify these functions with the kernel on the sphere, recall the transfer map \eqref{eq:transfer} and write
$y=(\cos(r/2)\theta,\sin(r/2)\omega)$, with
$\theta,\omega\in\s^{k-1}$. Direct differentiation yields
\begin{equation}\label{eq:maps}
 \begin{aligned}
 \mathcal T[Z_0]&=-h_{01}\omega_k,\\
 \mathcal T[Z_i]&=-h_{10}\theta_i+2w'\theta_i\omega_k,
                    &&1\le i\le k,\\
 \mathcal T[Z_{k+a}]&=-h_{01}\omega_a,
                    &&1\le a\le k-1,\\
 \mathcal T[Z_{n+i}]&=-h_{10}\theta_i-2w'\theta_i\omega_k,
                    &&1\le i\le k,\\
 \mathcal T[Z_{n+k+(i-1)(k-1)+a}]&=2w'\theta_i\omega_a,
                    &&\substack{1\le i\le k,\\1\le a\le k-1.}
 \end{aligned}
\end{equation}
In particular, for $1\le i\le k$,
\[
 \mathcal T[Z_i+Z_{n+i}]=-2h_{10}\theta_i,
 \qquad
 \mathcal T[Z_i-Z_{n+i}]=4w'\theta_i\omega_k.
\]
These identities show that the images of the functions in
\eqref{eq:euclidean-basis} span precisely the $k(k+2)$-dimensional kernel stated in Theorem~\ref{thm:main}. Observe that, for each
$1\le i\le k$, we may replace the pair
$Z_i,Z_{n+i}$ by $Z_i+Z_{n+i},Z_i-Z_{n+i}$,
since this change is invertible.
By \eqref{eq:maps}, the images under $\mathcal T$
of the resulting functions are nonzero multiples of
\[
\{h_{10}(r)\theta_i\}_{i=1}^k,\qquad
\{h_{01}(r)\omega_a\}_{a=1}^k,\qquad
\{w'(r)\theta_i\omega_a\}_{i,a=1}^k.
\]
These $k(k+2)$ functions are linearly independent,
because they are not identically
zero and form an orthogonal family.
Since $\mathcal T$ is invertible, the original
$k(k+2)$ functions in \eqref{eq:euclidean-basis}
are also linearly independent.

\begin{proof}[Proof of Theorem~\ref{thm:euclidean}]
The functions in \eqref{eq:euclidean-basis} are linearly independent
kernel elements. Conversely, \eqref{eq:transfer} transfers every
finite-energy kernel element in $\R^n$ to another kernel element in $\s^n$.
Theorem~\ref{thm:main} and \eqref{eq:maps} show that this image is a
linear combination of the images of the functions $Z_\alpha$.
The conclusion follows from the invertibility of the map defined in \eqref{eq:transfer}.
\end{proof}

\begin{proof}[Proof of Theorem~\ref{thm:B}]
Part (i) is Lemma~\ref{lem:unique-profile} and part (ii) is Proposition~\ref{prop:radial}. 
%
\end{proof}
\section*{Acknowledgements}
This work was initiated during the first author's visit to the Abdus Salam International Centre for Theoretical Physics (ICTP) in Trieste, Italy, supported by the ICTP–INdAM Collaborative Grants and Research in Pairs programme. The first author thanks ICTP for its support and warm hospitality.
 {The second author is a member of the INdAM research group GNAMPA.}

\section*{Declaration of generative AI}

During the preparation of this work the authors used ChatGPT (OpenAI) to assist with selected calculations and with drafting and revising the manuscript. After using this tool/service, the authors reviewed and edited the content as needed and take full responsibility for the content of the published article.

\bibliographystyle{elsarticle-harv}
\bibliography{biblio}

\end{document}